\documentclass[11pt,reqno,letterpaper]{amsart}
\usepackage[T1]{fontenc}
\usepackage{lmodern}
\usepackage{amsmath,amssymb,amsthm}
\usepackage[hmargin=1.12in,vmargin=1.07in]{geometry}
\usepackage[expansion=false]{microtype}
\usepackage[hidelinks]{hyperref}
\usepackage{color}
\hypersetup{pdftitle={The sharp logarithmic loss in vanishing order for real planar Schr\"odinger equations},
  pdfauthor={Igor Kukavica and Qi Xu},
  pdfsubject={Quantitative unique continuation for real potentials},
  pdfkeywords={vanishing order, unique continuation, real potential, closed surface, Landis conjecture}}

\newtheorem{theorem}{Theorem}
\newtheorem{corollary}[theorem]{Corollary}
\newtheorem{lemma}[theorem]{Lemma}
\numberwithin{equation}{section}
\DeclareMathOperator{\ord}{ord}
\DeclareMathOperator{\supp}{supp}
\newcommand{\R}{\mathbb R}
\newcommand{\C}{\mathbb C}
\newcommand{\D}{\mathbb D}
\newcommand{\dbar}{\bar\partial}
\renewcommand{\Re}{\operatorname{\mathbb{R}e}}
\newcommand{\comma}{,\qquad}
\newcommand{\andand}{\quad\text{and}\quad}

  \chardef\forshowkeys=0
  \chardef\showllabel=0
  \chardef\refcheck=0
  \chardef\sketches=0
  \chardef\figures=1

\ifnum\forshowkeys=1

  \usepackage[notref,notcite,color]{showkeys}
\fi

\definecolor{mygray}{rgb}{.6,.6,.6}
\definecolor{myblue}{rgb}{0.9, 0, 1}
\definecolor{colorforkeys}{rgb}{1.0,0.0,0.0}

\def\andand{\text{\indeq and\indeq}}

\def\supp{\mathop{\rm supp}\nolimits}
\def\indeq{\quad{}}

\definecolor{colorigor}{rgb}{.5, 0.2, 0.8}

\definecolor{colorgggg}{rgb}{0.1,0.5,0.3}
\definecolor{colorllll}{rgb}{0.0,0.7,0.0}
\definecolor{colorhhhh}{rgb}{0.3,0.75,0.4}
\definecolor{colorpppp}{rgb}{0.7,0.0,0.2}
\definecolor{coloroooo}{rgb}{0.45,0.0,0.0}
\definecolor{colorqqqq}{rgb}{0.1,0.7,0}

\definecolor{coloraaaa}{rgb}{0.6,0.6,0.6}

\def\comma{ {\rm ,\qquad{}} }

\newcommand{\eqnb}{\begin{equation}}
\newcommand{\eqne}{\end{equation}}

\title[Sharp vanishing orders for real planar Schr\"odinger equations]
{The sharp logarithmic loss in vanishing order\\ for real planar Schr\"odinger equations}
\author{Igor Kukavica}
\author{Qi Xu}
\address{Department of Mathematics, University of Southern California, Los Angeles, CA 90089}
\email[Igor Kukavica]{\href{mailto:kukavica@usc.edu}{kukavica@usc.edu}}
\email[Qi Xu]{\href{mailto:xuqi@usc.edu}{xuqi@usc.edu}}
\date{}

\begin{document}
\begin{abstract}
We determine, up to constants, the largest possible order of vanishing of real solutions of $\Delta u=Vu$ in dimension two in terms of the size of the real potential~$V$. For each sufficiently large integer $k$, we construct a real potential $V\in C_c^\infty(B_1)$ with $\|V\|_\infty$ comparable to $k^2/\log k$ and a smooth real solution $u$ on $\R^2$ that vanishes to order $k$ at the origin, satisfies $\sup_{B_{1/2}}|u|=\|u\|_\infty=1$, and equals a positive constant outside~$B_{7/8}$. Consequently, on every smooth connected closed Riemannian surface, the largest order of vanishing at a prescribed point for real potentials of norm at most $M$ is comparable to $\sqrt{M\log M}$ for large~$M$, uniformly in the point. This gives a negative answer, in dimension two, to the square-root questions of Kukavica and Kenig. With the same building block, we also construct a nonzero real solution on $\R^2$ with $\|V\|_\infty\le1$ and $|u(x)|\le\exp(-c|x|\sqrt{\log|x|})$ for large~$|x|$, so that the Landis conjecture fails in its exponential form in the plane even for real potentials. Both constructions match the estimates of Logunov, Malinnikova, Nadirashvili, and Nazarov, which are therefore sharp.
\end{abstract}
\maketitle

\section{Introduction}\label{EQ01}
We consider real solutions of the stationary Schr\"odinger equation
\begin{align}
 \Delta u=Vu\label{EQ02}
\end{align}
with a bounded real potential~$V$, and we address two questions of quantitative unique continuation in dimension two. The first is how large the order of vanishing of a nonzero solution can be in terms of~$\|V\|_\infty$. The second is how fast a nonzero solution on the whole plane can decay at infinity when $|V|\le1$.

The first question originates in the theorem of Donnelly and Fefferman \cite{DF}, who proved that an eigenfunction of $-\Delta_g$ with eigenvalue $\lambda$ on a smooth connected closed Riemannian manifold vanishes to order at most $C_g(1+\sqrt{\lambda})$ at every point. The power of $\lambda$ is optimal, since on the unit sphere $\mathbb S^2$ the spherical harmonic $\Re\big((x_1+i x_2)^k\big)$ has eigenvalue $k(k+1)$ and vanishes to order $k$ at the north pole. In 1998, the first named author asked in~\cite{Kukavica} whether the vanishing order on a closed manifold is at most $C_g(1+\sqrt{\|V\|_\infty})$ when $\lambda$ is replaced by a bounded real potential~$V$. In 2005, Kenig \cite[Question~2]{Kenig} independently posed a local version of this question. Suppose that $u$ and $V$ are real and satisfy $\Delta u+Vu=0$ in the ball $\{|x|<10\}\subset\R^n$, with $\|V\|_\infty\le M$, $\|u\|_\infty\le C_0$, and $\sup_{|x|<1}|u|\ge1$. Is it true that
\begin{align}
 \sup_{|x|<r}|u|\ge a_1r^{a_2\sqrt{M}}\label{EQ03}
\end{align}
for all sufficiently small $r>0$, where $a_1,a_2>0$ depend only on $n$ and~$C_0$?

Without the assumption that $u$ and $V$ are real, the answer is known. Kenig observed in~\cite{Kenig} that the Carleman estimate of Bourgain and Kenig~\cite{BK} gives \eqref{EQ03} with $M^{2/3}$ in place of~$\sqrt M$, and that the construction of Meshkov~\cite{Meshkov} shows that the power $2/3$ cannot be lowered for complex-valued $u$ and~$V$; see also~\cite{BakriComplex,CKW1}. For real potentials, Feng, Gan, Jiang, and Lin~\cite{FGJL} have recently shown that the power $2/3$ is optimal in all dimensions $n\ge3$. When $n=2$, Logunov, Malinnikova, Nadirashvili, and Nazarov~\cite{LMNN} proved an estimate, stated in \eqref{EQ08} below, which implies \eqref{EQ03} with $\sqrt{M\log M}$ in place of $\sqrt M$ for real solutions in dimension two. Thus in the plane the square-root bound can fail at most by a logarithmic factor, and it remained open whether this factor is necessary. We show that $\sqrt{M\log M}$ is the sharp rate in dimension two.

From here on, we identify $\R^2$ with $\C$ and write $B_r(a)=\{z\in\C:|z-a|<r\}$, $B_r=B_r(0)$, and $\D=B_1$. We denote $\Delta=\partial_x^2+\partial_y^2$ and use the same sign for~$\Delta_g$. In Euclidean or geodesic balls the vanishing order is
$$
 \ord_p u=\sup\biggl\{\alpha\ge0:\sup_{B_r(p)}|u|=O(r^\alpha)\text{ as }r\downarrow0\biggr\}.
$$
For nonnegative quantities we write $A\lesssim B$ if $A\le CB$ with a constant uniform in the parameters under consideration, and $A\sim B$ if both $A\lesssim B$ and $B\lesssim A$. Our main result is the following.

\begin{theorem}[The planar construction]\label{EQ04}
There are absolute constants $c,C>0$ such that for every sufficiently large integer $k$ there exist a real potential $V_k\in C_c^\infty(\D)$ and a real solution $u_k\in C^\infty(\R^2)$ satisfying
\begin{align}
 \Delta u_k=V_k u_k\quad\text{on }\R^2 \andand \ord_0 u_k=k.\label{EQ05}
\end{align}
Moreover,
\begin{align}
 \sup_{B_{1/2}}|u_k|=\|u_k\|_{L^\infty(\R^2)}=1 \andand
 c\frac{k^2}{\log k}\le\|V_k\|_\infty\le C\frac{k^2}{\log k}.\label{EQ06}
\end{align}
The solution is a positive constant on $\{|z|\ge7/8\}$ and satisfies
\begin{align}
 \sup_{B_r}|u_k|\le(4r)^k \comma 0<r\le1/4.\label{EQ07}
\end{align}
\end{theorem}

Theorem~\ref{EQ04} is optimal in view of the following estimate. For a bounded nonzero real solution of $\Delta u=Vu$ in $B_1$ with real $V$ and $\|V\|_\infty\le M$, where $M\ge2$, Logunov, Malinnikova, Nadirashvili, and Nazarov \cite[Corollary~2.3]{LMNN} proved
\begin{align}
 \sup_{B_r}|u|\ge r^{C(\sqrt{M\log M}+\gamma)}\sup_{B_1}|u|
 \comma 0<r\le\frac14,\label{EQ08}
\end{align}
where $C$ is an absolute constant and $\gamma=\log\big(\sup_{B_1}|u|/\sup_{B_{1/2}}|u|\big)$. In particular,
$$
 \ord_0 u\le C\big(\sqrt{M\log M}+\gamma\big).
$$
For the solutions in Theorem~\ref{EQ04}, we have $\gamma=0$. With $M=\|V_k\|_\infty$, the estimate~\eqref{EQ06} gives $\log M=2\log k-\log\log k+O(1)$, so $k\sim\sqrt{M\log M}$. Thus \eqref{EQ08} is sharp up to the value of the constant~$C$. Conversely, \eqref{EQ08} shows that the potentials in Theorem~\ref{EQ04} have the smallest possible size. Namely, a real potential on $B_1$ that admits a bounded real solution with $\gamma=0$ and vanishing order $k$ at the origin has norm at least $ck^2/\log k$ for large~$k$.

The solutions $u_k$ also settle Kenig's question in the plane. They satisfy its hypotheses with $C_0=1$, and
$$
 \sup_{B_1}|u_k|=\|u_k\|_{L^\infty(B_{10})}=1
 \andand
 \frac{\ord_0 u_k}{\sqrt{\|V_k\|_\infty}}\ge c\sqrt{\log k}\longrightarrow\infty.
$$
Indeed, for $M=\|V_k\|_\infty$ and any fixed $A>0$, once $k>A\sqrt M$, the estimate~\eqref{EQ07} gives $\sup_{B_r}|u_k|/r^{A\sqrt M}\to0$ as $r\downarrow0$, so that \eqref{EQ03} fails for every choice of $a_1$ and~$a_2$. On the other hand, if $u$ satisfies the hypotheses of Kenig's question in the plane, then \eqref{EQ08} applied to $x\mapsto u(2x)$ gives \eqref{EQ03} with $\sqrt{M\log M}$ in place of~$\sqrt M$. Hence, the optimal exponent in the local formulation in dimension two is $\sqrt{M\log M}$.

Since $u_k$ is constant near~$\partial\D$, the planar example can be inserted into an isothermal chart about any point of a closed surface.

\begin{theorem}[Closed surfaces]\label{EQ09}
Let $(X,g)$ be a smooth connected closed Riemannian surface. There exist $C_g,r_g>0$ and $k_g\in\mathbb N$ such that for every $p\in X$ and every integer $k\ge k_g$, one can find real $U_k,Q_k\in C^\infty(X)$ solving $\Delta_g U_k+Q_kU_k=0$ with
\begin{align}
  \ord_p U_k=k \andand
  \|Q_k\|_\infty\le C_g\frac{k^2}{\log k}.\label{EQ10}
\end{align}
The solution is normalized by $\|U_k\|_\infty=1$ and satisfies
\begin{align}
  \sup_{B_r(p)}|U_k|\le r^{k/2} \comma 0<r\le r_g.\label{EQ11}
\end{align}
It is a positive constant outside a coordinate disk about $p$, and $Q_k$ is compactly supported in that disk.
\end{theorem}

\begin{corollary}[Sharp order]\label{EQ12}
On each such surface there are $c_g,C_g>0$ and $M_g>2$ such that for every $p\in X$ and every $M\ge M_g$,
\begin{align}
  c_g\sqrt{M\log M}\le\sup_{U,Q}\ord_p U\le C_g\sqrt{M\log M},\label{EQ13}
\end{align}
where the supremum is over real $U,Q\in C^\infty(X)$ with $\Delta_g U+QU=0$, $\|U\|_\infty=1$, and $\|Q\|_\infty\le M$.
\end{corollary}

The upper bound in~\eqref{EQ13} is due to Logunov, Malinnikova, Nadirashvili, and Nazarov~\cite[Corollary~2.5]{LMNN}, while the lower bound follows from Theorem~\ref{EQ09}. Corollary~\ref{EQ12} answers the question of~\cite{Kukavica} in dimension two. On every closed surface, the bound $C_g(1+\sqrt{\|V\|_\infty})$ fails, but only by the factor~$\sqrt{\log\|V\|_\infty}$.

We now turn to the decay at infinity. According to~\cite{KL}, Landis conjectured in the late 1960s that a solution of \eqref{EQ02} in $\R^n$ with $|V|\le1$ that satisfies $|u(x)|\le\exp(-C|x|)$ for a sufficiently large constant $C$ vanishes identically. A weaker form of the conjecture asserts the same under the assumption $|u(x)|\le\exp(-|x|^{1+\varepsilon})$ for some $\varepsilon>0$. Meshkov~\cite{Meshkov} disproved both forms for complex potentials by constructing a nonzero solution in the plane with $|u(x)|\le\exp(-c|x|^{4/3})$, and he proved that the exponent $4/3$ cannot be increased. For real potentials, Frank and Ivanisvili~\cite{FI} recently obtained the same rate of decay in all dimensions $n\ge3$. In the plane, however, Logunov, Malinnikova, Nadirashvili, and Nazarov \cite[Theorem~1.1]{LMNN} proved the weak form of the conjecture for real potentials. Namely, if $u$ and $V$ are real, $|V|\le1$, and
\begin{align}
 |u(x)|\le\exp\big(-C|x|\sqrt{\log|x|}\big)\comma |x|>2,\label{EQ14}
\end{align}
with a sufficiently large absolute constant~$C$, then $u\equiv0$. Whether the factor $\sqrt{\log|x|}$ can be removed, that is, whether the conjecture holds in its original exponential form for real planar potentials, remained open. The next theorem, which is based on the same one-pole construction as Theorem~\ref{EQ04}, shows that it cannot.

\begin{theorem}[Decay at infinity]\label{EQ15}
There exist real $u,V\in C^\infty(\R^2)$ satisfying $\Delta u=Vu$ on $\R^2$, with $u(0)=\|u\|_\infty=1$ and $\|V\|_\infty\le1$, such that for some $c>0$ and $r_0>2$,
\begin{align}
 \sup_{|x|=r}|u(x)|\le\exp\big(-cr\sqrt{\log r}\big)\comma r\ge r_0.\label{EQ16}
\end{align}
In particular, $e^{A|x|}|u(x)|\to0$ as $|x|\to\infty$ for every $A>0$.
\end{theorem}

Thus the exponential form of the Landis conjecture fails in the plane even for real potentials. Moreover, a nonzero real solution with $|V|\le1$ can decay as fast as $\exp(-c|x|\sqrt{\log|x|})$ when $c>0$ is small, while \eqref{EQ14} with a large constant $C$ forces $u\equiv0$. This should be compared with the rate $\exp(-c|x|^{4/3})$, which is optimal for complex potentials in the plane~\cite{Meshkov} and for real potentials in dimensions $n\ge3$~\cite{FI}.

Theorem~\ref{EQ04} shows that the local version of the result in~\cite{LMNN} is sharp as well. By \cite[Theorem~1.2]{LMNN}, a real solution of $\Delta v=Wv$ on $\R^2$ with real $W$ and $|W|\le1$, for which the maximum $\|v\|_\infty=1$ of $|v|$ is attained at a point $y$ with $R/4\le|y|\le R/2$, satisfies $\sup_{B_1}|v|\ge\exp(-CR\log^{3/2}R)$ for large~$R$. On the other hand, let $k=\lfloor\varepsilon R\sqrt{\log R}\rfloor$, where $\varepsilon>0$ is a small absolute constant, and set $v(x)=u_k(x/R)$. By \eqref{EQ06} and \eqref{EQ07}, the function $v$ is such a solution, and $\sup_{B_1}|v|\le(4/R)^k\le\exp(-cR\log^{3/2}R)$.

\subsection*{Related work}
Unique continuation for \eqref{EQ02} goes back to Carleman's planar uniqueness theorem \cite{Carleman} and Aronszajn's strong unique continuation theorem~\cite{Aronszajn}. Jerison and Kenig \cite{JK} established unique continuation for Schr\"odinger operators with rough potentials, and Garofalo and Lin \cite{GL} obtained doubling estimates by a frequency argument. The square-root bound for the order of vanishing does hold with derivative control on the potential, as Bakri~\cite{BakriC1} proved on closed manifolds and Zhu \cite{Zhu} proved locally. In the plane, Kenig, Silvestre, and Wang~\cite{KSW} obtained a square-root bound in terms of $\|V\|_\infty$ for $\Delta u=Vu$ with $V\ge0$ under a two-scale normalization. The example in Theorem~\ref{EQ04} shows what is lost when the sign restriction is removed. For complex potentials, Bakri \cite{BakriComplex} attained order $M^{2/3}$ on $\mathbb S^2$, and Camliyurt, Kukavica, and Wang~\cite{CKW1} attained the same order in planar examples with periodic or Dirichlet boundary conditions; see also~\cite{CKW2}. For real potentials, Feng, Gan, Jiang, and Lin \cite[Theorem~1.1 and Remark~1.2]{FGJL} constructed Dirichlet solutions on a three-dimensional ball with vanishing order $k$ and $\|V_k\|_\infty\le Ck^{3/2}$. Their construction extends to product domains in all higher dimensions. For bounds on the order of vanishing when the lower order coefficients belong to Lebesgue classes, see Davey and Zhu~\cite{DZ} and Davey~\cite{Davey}.

Concerning the decay at infinity, Bourgain and Kenig \cite{BK} obtained a quantitative lower bound in their work on Anderson localization, which matches the example of Meshkov~\cite{Meshkov} up to a logarithm in the exponent. For real planar potentials, Kenig, Silvestre, and Wang~\cite{KSW} treated the case $V\ge0$, and Davey, Kenig, and Wang \cite{DKW} treated the case of an exponentially decaying negative part, before the sign restriction was removed in~\cite{LMNN}. Le Balc'h and Souza \cite{LBS} extended planar quantitative estimates to equations with bounded real drift terms. In the periodic cylinder $\R\times(0,2\pi)^d$ with the axial variable~$t$, Filonov and Krymskii~\cite{FK} showed that the fastest decay of a nonzero real solution with a bounded real potential is of the form $e^{-c|t|}$ for $d=1,2$ and $e^{-c|t|^{4/3}}$ for $d\ge3$. We refer to the survey~\cite{FSR} for further results and references on the Landis conjecture.

\subsection*{Outline of the construction}
In Section~\ref{EQ17}, we introduce a rational function $R$ with a zero of order $k$ at the origin and $k$ simple poles in the annulus $\{1/4<|z|<1/2\}$, separated by distances of at least $ck^{-1/2}$, such that $R-1$ is exponentially small for $|z|\ge3/4$. Near each pole, $R$ is a simple pole multiplied by an exponential whose frequency may be comparable to~$k$. The main step is Lemma~\ref{EQ31}, which replaces $\Re R$ near a pole by a smooth solution of $\Delta u=Vu$ with a potential of size $k^2/\log k$. We first lower the frequency to at most $k/2\sqrt{\log k}$ and remove the pole at the scale $\sqrt{\log k}/k$. We then restore the frequency across an annulus of logarithmic width comparable to~$\log k$, where a real amplitude correction makes the Laplacian a smooth multiple of the solution even at its zeros. Section~\ref{EQ54} gives the normalization and the constant exterior value, and Section~\ref{EQ57} transfers the example to a closed surface. In Section~\ref{EQ59}, we apply the same lemma to an infinite product with $j$ poles on the circle of radius $j^2/\sqrt{\log j}$, which gives Theorem~\ref{EQ15}.

\section{The regular-polygon product}\label{EQ17}
In the paper, constants are absolute unless a dependence is indicated. For a sufficiently large integer $k$, let $N=\lfloor\sqrt{k}\rfloor$ and $m=\lfloor k/N\rfloor$. For $1\le j\le N$, set
$$
 \rho_j=\frac{1}{4}+\frac{j}{4(N+1)}
 \andand n_j=\left\lfloor\frac{k+j-1}{N}\right\rfloor,
$$
and consider the rational function
\begin{align}
 R(z)=\prod_{j=1}^N\frac{z^{n_j}}{z^{n_j}-\rho_j^{n_j}}.\label{EQ18}
\end{align}
Note that its poles on $|z|=\rho_j$ are the vertices of a regular
$n_j$-gon. Write $k=Nm+r$, where $0\le r<N$. Then
\begin{align*}
n_j
=\left\lfloor\frac{Nm+r+j-1}{N}\right\rfloor
=\begin{cases}
m, & 1\le j\le N-r,\\
m+1, & N-r<j\le N.
\end{cases}
\end{align*}
Since $N^2\le k\le N^2+2N$,
$$
 \sum_{j=1}^N n_j=(N-r)m+r(m+1)=k \andand N\le n_j\le N+2.
$$
All poles are simple.
For distinct poles $a,a'$ on circles $i,j$, respectively,
\begin{align*}
|a-a'|
&\ge|\rho_i-\rho_j|=\frac{|i-j|}{4(N+1)}
\ge\rho_2-\rho_1,
&&i\ne j,\\
|a-a'|
&\ge2\rho_j\sin\frac{\pi}{n_j}
\ge\frac{4\rho_j}{n_j}
\ge\frac{1}{N+2}
>\rho_2-\rho_1,
&&i=j.
\end{align*}
Equality in the first case is attained at $a=\rho_1$ and
$a'=\rho_2$. Hence,
\begin{align}
\min_{a\ne a'}|a-a'|
=\rho_2-\rho_1
=\frac{1}{4(N+1)}
=\frac{1+O(k^{-1/2})}{4\sqrt{k}}.
\label{EQ19}
\end{align}

We first estimate $R$ near the origin, on the outer
annulus, and near each pole.

\smallskip\noindent\textit{At the origin.}
Set $c_k=(-1)^N\prod_j\rho_j^{-n_j}$. For $|z|<\rho_1$, we have $R(z)=c_kz^k\prod_j(1-(z/\rho_j)^{n_j})^{-1}$. Expanding the factors at zero gives
\begin{align}
 R(z)=c_k z^k\biggl(1+z^m\sum_{n_j=m}\rho_j^{-m}+O_k(|z|^{m+1})\biggr)
 \indeq \text{as~}z\to0.\label{EQ20}
\end{align}
Here, $c_k\in\R\setminus\{0\}$ and $\sum_{n_j=m}\rho_j^{-m}
=\sum_{j=1}^{N-r}\rho_j^{-m}>0$. Bernoulli's inequality and $n_j\ge N$ give
$$
 (4\rho_j)^{n_j}-1
 =\left(1+\frac{j}{N+1}\right)^{n_j}-1\ge\frac{Nj}{N+1}\ge\frac{j}{2}.
$$
For $|z|\le1/4$ and $N\ge4$, the reverse triangle inequality in \eqref{EQ18} yields
\begin{align}
 |R(z)|\le\prod_{j=1}^N\frac{(4|z|)^{n_j}}{(4\rho_j)^{n_j}-1}
 \le\frac{2^N}{N!}(4|z|)^k\le(4|z|)^k.\label{EQ21}
\end{align}

\smallskip\noindent\textit{On the outer annulus.}
For $|z|\ge3/4$, let $t_j(z)=(\rho_j/z)^{n_j}$. Since $\rho_j<1/2$ and $n_j\ge N$, we have $|t_j|\le(2/3)^N\le1/2$ for large~$k$. From $R=\prod_j(1-t_j)^{-1}$, the branch of $\log R$ that vanishes at infinity satisfies
$$
 |\log R(z)|=\left|-\sum_{j=1}^N\log(1-t_j)\right|\le\sum_{j=1}^N\frac{|t_j|}{1-|t_j|}\le2N(2/3)^N.
$$
Using $e^w-1=w\int_0^1e^{tw}\,\text{d}t$ with $w=\log R(z)$ gives
\begin{align}
 |R(z)-1|\le e^{|\log R(z)|}|\log R(z)|\lesssim N(2/3)^N
    \comma |z|\ge3/4.
\label{EQ22}
\end{align}
Since $t_j'=-n_jt_j/z$, $\sum_j n_j=k$, and $|R|\le e^{|\log R|}$, logarithmic differentiation yields
\begin{align}
 |R'(z)|=|R(z)|\left|\sum_{j=1}^N\frac{-n_jt_j}{z(1-t_j)}\right|\le\frac{2|R(z)|}{|z|}\sum_{j=1}^N n_j|t_j|\lesssim k(2/3)^N
    \comma |z|\ge3/4.
    \label{EQ23}
\end{align}

\smallskip\noindent\textit{Near a pole.}
Let
$d=(\rho_2-\rho_1)/16=1/[64(N+1)]$, so that $\rho_1-1/4=1/2-\rho_N=16d>8d$.
By \eqref{EQ19}, the closed disks $\overline{B_{4d}(a)}$
about the poles are pairwise disjoint and contained in
$\{1/4<|z|<1/2\}$.

Fix a pole $a$ on the $j$th circle, so $|a|=\rho_j$.
With $c_a=\operatorname{Res}_a R$, the function
$(z-a)R(z)/c_a$ is holomorphic and nonvanishing near
$\overline{B_{4d}(a)}$ and equals one at~$a$.
Choose its holomorphic logarithm $p$ on a neighborhood of this closed disk with $p(a)=0$.
We claim that
\begin{align}
p'(z)=\frac1z\sum_{i>j}n_i+O(N),
\qquad |z-a|\le4d.
\label{EQ24}
\end{align}
Writing $R_i(z)=z^{n_i}/z^{n_i}-\rho_i^{n_i}$ for the $i$-th factor in~\eqref{EQ18}, we have 
$$
p'(z)
=\frac1{z-a}+\sum_{i=1}^N\frac{R_i'(z)}{R_i(z)}
=\underbrace{\left(\frac1{z-a}
+\frac{R_j'(z)}{R_j(z)}\right)}_{\text{regularized $j$th factor}}
+\underbrace{\sum_{i\ne j}\frac{R_i'(z)}{R_i(z)}}_{\text{remaining factors}}
$$
for $0<|z-a|\le4d$.

\smallskip\noindent\textit{The factors with $i\ne j$.}
For $i\ne j$, the logarithmic derivative of the
$i$th factor of $R$ is
$$
\frac{R_i'(z)}{R_i(z)}
=\frac{\text{d}}{\text{d}z}\log R_i(z)
=\frac{n_i}{z}
-\frac{n_i z^{n_i-1}}{z^{n_i}-\rho_i^{n_i}}
=\frac{n_i}{z\left(1-(z/\rho_i)^{n_i}\right)}.
$$
For $|z-a|\le4d$ and $i\ne j$, the separation of
the circles gives
$$
\bigl||z|-\rho_i\bigr|\ge|\rho_j-\rho_i|-|z-a|
\ge\frac{|i-j|}{4(N+1)}-4d
\ge\frac{3|i-j|}{16(N+1)}.
$$
Since $n_i\ge N$ and $\max\{|z|,\rho_i\}<1/2$, we obtain for $N\ge2$
$$
n_i\left|\log\frac{|z|}{\rho_i}\right|=n_i|\log|z|-\log\rho_i|
\ge\frac{n_i\big||z|-\rho_i\big|}
{\max\{|z|,\rho_i\}}
\ge\frac{3N|i-j|}{8(N+1)}
\ge\frac{|i-j|}{4}.
$$
In particular, $|z|<\rho_i$ if $i>j$, while
$|z|>\rho_i$ if $i<j$.
Consequently, for $t=(z/\rho_i)^{n_i}$,
$$
|t|\le e^{-|i-j|/4}
    \comma i>j
$$
and
$$
|t|\ge e^{|i-j|/4}
    \comma i<j
    .
$$
Since $|i-j|\ge1$, the two cases give
\begin{align*}
\left|\frac1{1-t}-\mathbf1_{\{i>j\}}\right|
&\le
\begin{cases}
\dfrac{|t|}{1-|t|},&i>j,\\
\dfrac1{|t|-1},&i<j,
\end{cases}\le\frac1{e^{|i-j|/4}-1}
\leq \frac{1}{1-e^{-1/4}} e^{-|i-j|/4}.
\end{align*}
Therefore, summing over $i\ne j$ and using
$n_i\le N+2$ and $|z|\ge1/4$, we obtain
\begin{align}
\begin{split}
\label{EQ25}
\left|
\sum_{i\ne j}\frac{R_i'(z)}{R_i(z)}
-\frac1z\sum_{i>j}n_i
\right|
&\le\frac1{|z|}\sum_{i\ne j}n_i
\left|\frac1{1-(z/\rho_i)^{n_i}}
-\mathbf1_{\{i>j\}}\right|\\
&\lesssim N\sum_{\ell=1}^{\infty}
e^{-\ell/4}
\lesssim N,
\qquad |z-a|\le4d.
\end{split}
\end{align}

\smallskip\noindent\textit{The factor with $i=j$.}
For $z\ne a$, we have
\begin{align}
\label{EQ26}
\frac1{z-a}+\frac{R_j'(z)}{R_j(z)}
=\frac{\text{d}}{\text{d}z}\log\bigl((z-a)R_j(z)\bigr).
\end{align}
For $n=n_j$, we have $a^n=\rho_j^n$. By the fundamental
theorem of calculus,
$$
z^n-a^n
=n(z-a)\int_0^1(a+t(z-a))^{n-1}\,\text{d}t=n(z-a)a^{n-1}I(z),
$$
where
$$
I(z)=\int_0^1
\left(1+t\frac{z-a}{a}\right)^{n-1}\,\text{d}t.
$$
Therefore,
\begin{align}
\label{EQ27}
(z-a)R_j(z)
=(z-a)\frac{z^n}{z^n-a^n}
=\frac{a}{n}\frac{(z/a)^n}{I(z)}.
\end{align}
Since $n\le N+2$, $d=1/[64(N+1)]$, and $|a|=\rho_j$, we have
$$
 \frac{4nd}{|a|}=\frac{n}{4(N+1+j)}\le\frac14.
$$
For any $w\in\mathbb C$, the binomial theorem gives
$$
|(1+w)^{n-1}-1|
\le\sum_{\ell=1}^{n-1}
\binom{n-1}{\ell}|w|^\ell
=(1+|w|)^{n-1}-1
\le e^{(n-1)|w|}-1.
$$
Applying this inequality to $I(z)$, 
we obtain
$$
|I(z)-1|
\le\int_0^1
\left(e^{(n-1)t|z-a|/|a|}-1\right)\,\text{d}t\le e^{4nd/|a|}-1
\le e^{1/4}-1<\frac12
$$
for $|z-a|\le4d$.
Differentiating under the integral also gives
$$
|I'(z)|
=\left|\frac{n-1}{a}\int_0^1
t\left(1+t\frac{z-a}{a}\right)^{n-2}\,\text{d}t\right|
\le\frac{n-1}{|a|}
\int_0^1t\,e^{4nd/|a|}\,\text{d}t
\lesssim n.
$$
Since $|I(z)|\ge1/2$ and $|z|\ge1/4$, \eqref{EQ26} and \eqref{EQ27} give
\begin{align}
 \left|\frac1{z-a}+\frac{R_j'(z)}{R_j(z)}\right|
 =\left|\frac nz-\frac{I'(z)}{I(z)}\right|
 \le\frac n{|z|}+\frac{|I'(z)|}{|I(z)|}\lesssim n\lesssim N.\label{EQ28}
\end{align}
Combining \eqref{EQ25} and \eqref{EQ28} proves \eqref{EQ24} for $z\ne a$. The expression $n/z-I'(z)/I(z)$ extends holomorphically across $a$, so the estimate holds there as well.

Set $b=p'(a)$ and $h(z)=p(z)-b(z-a)$. Then $h(a)=h'(a)=0$, and \eqref{EQ24} gives
$$
 |b|\le\frac1{|a|}\sum_{i>j}n_i+CN\lesssim k.
$$
Applying \eqref{EQ24} at both $z$ and $a$, we obtain
\begin{align*}
|h'(z)|
&
=|p'(z)-p'(a)|
\le\left|\frac1z-\frac1a\right|
\sum_{i>j}n_i+CN
\\&
\le\frac{|z-a|}{|az|}\sum_{i>j}n_i+CN
\lesssim kd+N\lesssim N,
\qquad |z-a|\le4d.
\end{align*}
Since $h(a)=0$, integration along the segment from
$a$ to $z$ yields
$$
\|h\|_{L^\infty(B_{2d}(a))}
\le2d\|h'\|_{L^\infty(B_{2d}(a))}
\lesssim Nd\lesssim1.
$$
Thus for $0<|z-a|<2d$,
\begin{align}
 R(z)=\frac{c_a}{z-a}e^{b(z-a)+h(z)} \andand \frac{|b|}{k}+\|h\|_{L^\infty(B_{2d}(a))}\le C,\label{EQ29}
\end{align}
with a constant $C$ independent of $k$ and~$a$.

\section{One-pole regularization}\label{EQ30}
We smoothen $\Re R$ near each pole by introducing a real potential of size $k^2/\log k$.

\begin{lemma}[One-pole regularization]\label{EQ31}
Fix $A\ge1$ and $a\in\C$. For an integer $k\ge2$, set $d=1/[64(\lfloor\sqrt{k}\rfloor+1)]$. Suppose $h$ is holomorphic on $B_{2d}(a)$ with $h(a)=h'(a)=0$, and let
$$
 R(z)=\frac{c_a}{z-a}e^{b(z-a)+h(z)} \comma 0<|z-a|<2d,
$$
where $c_a\in\C\setminus\{0\}$, $b\in\C$, and $|b|/k+\|h\|_{L^\infty(B_{2d}(a))}\le A$. For $k\ge k_0(A)$ there are real $u\in C^\infty(\overline{B_d(a)})$ and $V\in C_c^\infty(B_d(a))$ satisfying $\Delta u=Vu$ in~$B_d(a)$. The pair agrees with $(\Re R,0)$ near $\partial B_d(a)$, and
$$
 c\frac{k^2}{\log k}\le\|V\|_\infty\le C_A\frac{k^2}{\log k}.
$$
Here, $c>0$ is absolute, and $C_A$ and $k_0(A)$ depend only on~$A$.
\end{lemma}

\begin{proof}
Translate the pole to $a=0$ and write $r=|z|$.
We first regularize the pole with a truncated linear coefficient
and then restore the original field away from the pole. Set
\begin{align}
 b_*=\frac{b}{\max\{1,\,2|b|\sqrt{\log k}/k\}}
 \andand F(z)=\frac{c_a}{z}e^{b_*z+h(z)}
 \comma 0<r\le4\sqrt{\log k}/k.
 \label{EQ32}
\end{align}
Since $d\sim k^{-1/2}$,
we have $4\sqrt{\log k}/k<d/2$ for large~$k$.
We proceed from the pole to $\partial B_d$ in four steps.
\begin{enumerate}
\item[(i)] \textbf{Core.} For $r<\sqrt{\log k}/k$, use the scaled reciprocal of $F$ to smooth the pole and obtain both bounds for the potential.
\item[(ii)] \textbf{Inner annulus.} Starting at $r=\sqrt{\log k}/2k$,
match the core pair to $(\Re F,0)$ and continue it harmonically to
$r=4\sqrt{\log k}/k$.
\item[(iii)] \textbf{Correction annulus.} Restore $b$ on
$4\sqrt{\log k}/k<r<d/2$ with a cutoff in~$\log r$. A real amplitude
correction gives a smooth potential across the zeros of~$\Re F$.
\item[(iv)] \textbf{Outer annulus.} Match $(\Re R,0)$ on $d/2<r\le d$.
\end{enumerate}

\smallskip\noindent\textbf{(i) Core.}
On $r<k^{-1}\sqrt{\log k}$ the scaled reciprocal of $F$ extends holomorphically
across the pole as
\begin{align}
 \xi(z)=\frac{kz}{\sqrt{\log k}}e^{-b_*z-h(z)}
 \andand F(z)=\frac{k c_a}{\sqrt{\log k}}\frac{\bar\xi}{|\xi|^2}
 \comma z\ne0.\label{EQ33}
\end{align}
To replace $|\xi|^{-2}$ by a smooth positive profile $q(|\xi|)$, first
compute the Laplacian. For $t=|w|>0$, we have
$\Delta_w(q(t)\bar w)=(q''(t)+t^{-1}q'(t))\bar w+2t^{-1}q'(t)\bar w$.
Since $\xi$ is holomorphic, $\Delta_z(v\circ\xi)=|\xi'|^2(\Delta_w v)\circ\xi$,
so
\begin{align}
 \Delta_z\big(q(|\xi|)\bar\xi\big)
 =|\xi'|^2\left(q''(|\xi|)+\frac{3q'(|\xi|)}{|\xi|}\right)\bar\xi
 \comma \xi\ne0.\label{EQ34}
\end{align}
Here, $q''+3q'/t=t^{-3}(t^3q')'$, and $t^3q'$ equals zero for a constant
profile and $-2$ for $q(t)=t^{-2}$. Fix once and for all a smooth nondecreasing
$\eta\colon\R\to[0,1]$ equal to zero on $(-\infty,0]$ and one on $[1,\infty)$.
Interpolate between these two values of $t^3q'$ by setting
\begin{align}
 q(t)=2\int_t^\infty\frac{\eta(8v-1)}{v^3}\,\text{d}v
 \comma t\ge0,
 \qquad\text{so that}\qquad t^3q'(t)=-2\eta(8t-1)
 \comma t>0.
 \label{EQ35}
\end{align}
Thus $q$ is positive, constant on $[0,1/8]$, and equal to $t^{-2}$ for
$t\ge1/4$. Let
\begin{align}
 u(z)=\frac{k}{\sqrt{\log k}}\,q(|\xi|)\Re(c_a\bar\xi).
 \label{EQ36}
\end{align}
By \eqref{EQ34} and \eqref{EQ35}, $\Delta u=Vu$ at every $\xi\ne0$, with
\begin{align}
 V=|\xi'|^2\frac{q''(|\xi|)+3q'(|\xi|)/|\xi|}{q(|\xi|)}
 =-\frac{16|\xi'|^2\eta'(8|\xi|-1)}{|\xi|^3q(|\xi|)}.
 \label{EQ37}
\end{align}
For $|\xi|<1/8$, \eqref{EQ36} is a constant multiple of
$\Re(c_a\bar\xi)$ and is harmonic. Taking $V(0)=0$ thus gives a smooth
pair throughout the core.

To estimate $|\xi'|$ in \eqref{EQ37}, apply Cauchy's estimate on disks
of radius~$d/2$. We have $|h'(z)|\le2A/d$ on $B_d$ and, by \eqref{EQ32},
\begin{align}
 \frac{\sqrt{\log k}}{k}\sup_{|z|\le\sqrt{\log k}/k}|b_*+h'(z)|
 \le\frac12+\frac{2A\sqrt{\log k}}{kd}\le\frac23\label{EQ38}
\end{align}
for sufficiently large~$k$. In particular $|z(b_*+h'(z))|\le2/3$ in the
core. Since $h(0)=0$, integration along the segment from $0$ to $z$ yields
\begin{align}
 |b_*z+h(z)|=\left|\int_0^1 z\big(b_*+h'(tz)\big)\,\text{d}t\right|
 \le\frac{2kr}{3\sqrt{\log k}}\le\frac23
 \comma r\le k^{-1}\sqrt{\log k}.
 \label{EQ39}
\end{align}
Differentiating \eqref{EQ33} and using \eqref{EQ38} and \eqref{EQ39} gives
\begin{align}
 \frac{e^{-2/3}}{3}
 \le\frac{\sqrt{\log k}}{k}|\xi'(z)|
 =|e^{-b_*z-h(z)}|\,|1-z(b_*+h'(z))|
 \le\frac53 e^{2/3}.\label{EQ40}
\end{align}
By \eqref{EQ37}, $V$ vanishes unless $1/8\le|\xi|\le1/4$. For
$1/8\le t\le1/4$, \eqref{EQ35} gives
$$
 16\le q(t)=16+2\int_t^{1/4}\frac{\eta(8v-1)}{v^3}\,\text{d}v\le t^{-2}.
$$
Hence, $|\xi|^3q(|\xi|)\ge16(1/8)^3=1/32$ on the support of $V$ in the
core. Equations \eqref{EQ37} and \eqref{EQ40} give
$|V|\le512\|\eta'\|_\infty|\xi'|^2\lesssim k^2/\log k$ there.

For the lower bound, observe from \eqref{EQ33} and \eqref{EQ39} that
\begin{align}
 |\xi(z)|=\frac{kr}{\sqrt{\log k}}e^{-\Re(b_*z+h(z))}
 \ge e^{-2/3}\frac{kr}{\sqrt{\log k}}
    \comma 
 r\le k^{-1}\sqrt{\log k}.\label{EQ41}
\end{align}
Along the positive real axis $|\xi|$ goes from zero at $z=0$ to at least
$e^{-2/3}/2>1/4$ at $z=\sqrt{\log k}/2k$.

Choose $\theta\in(0,1)$ with $\eta'(\theta)=1$ by the mean value theorem.
Continuity on this segment gives $z_0\in(0,\sqrt{\log k}/2k)$ with
$|\xi(z_0)|=(1+\theta)/8$. Since
$q(|\xi(z_0)|)\le|\xi(z_0)|^{-2}$, \eqref{EQ37} and \eqref{EQ40} yield
$$
 |V(z_0)|=\frac{16|\xi'(z_0)|^2}{|\xi(z_0)|^3q(|\xi(z_0)|)}
 \ge\frac{16|\xi'(z_0)|^2}{|\xi(z_0)|}
 \ge\frac{64e^{-4/3}}9\,\frac{k^2}{\log k}.
$$
Together with the upper bound, this gives
\begin{align}
 c\frac{k^2}{\log k}\le\sup_{kr<\frac12\sqrt{\log k}}|V|
 \le\sup_{kr<\sqrt{\log k}}|V|\le C\frac{k^2}{\log k}.
 \label{EQ42}
\end{align}

\smallskip\noindent\textbf{(ii) Inner annulus.}
On $\sqrt{\log k}/2k\le r\le\sqrt{\log k}/k$, the estimate~\eqref{EQ41} gives $|\xi|\ge e^{-2/3}/2>1/4$. Hence, $q(|\xi|)=|\xi|^{-2}$ there, and \eqref{EQ33}, \eqref{EQ36}, and \eqref{EQ37} give $(u,V)=(\Re F,0)$. Using the holomorphic field in \eqref{EQ32}, we may therefore set
$$
 (u,V)=(\Re F,0)
     \comma \frac{\sqrt{\log k}}{2k}\le r\le\frac{4\sqrt{\log k}}{k}.
$$

\smallskip\noindent\textbf{(iii) Correction annulus.}
If $b=b_*$, the field in \eqref{EQ32} is already $R$, and
$(\Re R,0)$ extends the pair from (ii) to~$\partial B_d$. For the rest of (iii), suppose $b\ne b_*$.
To extend $F$ from the inner annulus, set
\begin{align}
 \beta(r)=b_*+(b-b_*)\eta\left(
 \frac{\log(kr/4\sqrt{\log k})}{\log(kd/8\sqrt{\log k})}\right)
 \andand F(z)=\frac{c_a}{z}e^{z\beta(r)+h(z)}
 \label{EQ43}
\end{align}
for $0<r\le d$.
The cutoff equals 0 for $r\le4 k^{-1}\sqrt{\log k}$ and 1 for $r\ge d/2$.
Since $\beta$ lies on the segment from $b_*$ to $b$ on the same ray, on
$4\sqrt{\log k}/k\le r\le d$ we have
\begin{align}
 \frac{k}{2\sqrt{\log k}}=|b_*|\le|\beta|\le|b|\le Ak
 \andand r^{-1}\le\frac{|b_*|}{2}\le\frac{|\beta|}{2}.
 \label{EQ44}
\end{align}
We first show that $\nabla\Re F$ stays nonzero where $\beta$ varies.
With $\partial=(\partial_x-i\partial_y)/2$ and
$\dbar=(\partial_x+i\partial_y)/2$, differentiating \eqref{EQ43} gives
\begin{align}
 \frac{\partial F}{F}=\beta-z^{-1}+h'+z\partial\beta
 \andand \frac{\dbar F}{F}=z\dbar\beta.\label{EQ45}
\end{align}
For the radial function $\beta$, we have $|\partial\beta|=|\dbar\beta|=|\beta'(r)|/2$.
Since $d\sim k^{-1/2}$, the denominator of the cutoff argument in
\eqref{EQ43} is $\tfrac12\log k-\tfrac12\log\log k+O(1)$, so
$$
 r|\beta'(r)|\le
 \frac{|b-b_*|\,\|\eta'\|_\infty}{\log(kd/8\sqrt{\log k})}
 \lesssim_A\frac{k}{\log k}.
$$
Using $|\nabla\Re F|=|2\partial\Re F|=|\partial F+\overline{\dbar F}|$,
\eqref{EQ45}, and Cauchy's bound $|h'|\le2A/d$, we obtain
\begin{align}
 \begin{split}
 \frac{|\nabla\Re F|}{|F|}
 &\ge\frac{|\partial F|-|\dbar F|}{|F|}
 \ge|\beta|-r^{-1}-|h'|-r|\beta'(r)|\\
 &\ge\frac{|\beta|}{2}-C_A\left(d^{-1}+\frac{k}{\log k}\right)
 \ge\frac{|\beta|}{4}.
 \end{split}\label{EQ46}
\end{align}
Indeed, by \eqref{EQ44} the error divided by $|\beta|$ is at most
$C_A(\sqrt{\log k}/kd+1/\sqrt{\log k})$, which tends to zero.
By \eqref{EQ46}, $|\nabla\Re F|>0$ throughout the correction
annulus. Define
\begin{align}
 g=-\frac{(\Re F)\Delta\Re F}{2|\nabla\Re F|^2}
 \andand u=e^g\Re F.\label{EQ47}
\end{align}
The product rule gives
\begin{align*}
 e^{-g}\Delta u
 &=\Delta\Re F+2\nabla g\cdot\nabla\Re F
   +(\Re F)(\Delta g+|\nabla g|^2)\\
 &=(\Re F)\left[
   \Delta g+|\nabla g|^2
   -\nabla\Re F\cdot\nabla
   \left(\frac{\Delta\Re F}{|\nabla\Re F|^2}\right)
   \right].
\end{align*}
The identity $\Delta u=Vu$ thus holds also at the zeros of $\Re F$, with
\begin{align}
 V=\Delta g+|\nabla g|^2-\nabla\Re F\cdot\nabla
 \left(\frac{\Delta\Re F}{|\nabla\Re F|^2}\right).
 \label{EQ48}
\end{align}
At $r=4\sqrt{\log k}/k$, \eqref{EQ46} and continuity keep the denominator in \eqref{EQ47} and \eqref{EQ48} positive in a neighborhood of the circle. On its inner side $\beta=b_*$, so $F$ is holomorphic and $\Delta\Re F=0$. The same smooth formulas thus give $g=V=0$ on the inner side and join the pair in~(ii).

To bound the derivatives in \eqref{EQ48}, write $D^j$ for the array of
order-$j$ Cartesian derivatives, with the maximum norm. In a derivative
of the cutoff in \eqref{EQ43}, each factor $D^m\log r$ is
$O_m(r^{-m})$, and the orders of these factors sum to $j$. Hence,
\begin{align}
 r^j|D^j\beta|
 \lesssim_j |b-b_*|\sum_{\ell=1}^j
 \frac{\|\eta^{(\ell)}\|_\infty}{[\log(kd/8\sqrt{\log k})]^\ell}
 \lesssim_A\frac{k}{\log k},\qquad 1\le j\le4.\label{EQ49}
\end{align}
For $1\le j\le3$, Cauchy's estimate gives $|D^j h|\le j!A(2/d)^j$.
At most one derivative falls on $z$ in~$D^j(z\beta)$. Using \eqref{EQ49},
$d^{-1}\le r^{-1}\le|\beta|/2$, and $k/\log k\lesssim|\beta|$, we get
$$
 |D^j(z\beta+h)|\lesssim_A
 \mathbf1_{\{j=1\}}|\beta|+\frac{k}{\log k}r^{1-j}+d^{-j}
 \lesssim_A|\beta|^j.
$$
The chain rule and the preceding estimate give
$$
 |D^m(e^{z\beta+h})|
 \lesssim_A |e^{z\beta+h}|\,|\beta|^m,
 \qquad 0\le m\le3.
$$
Since $|D^i(z^{-1})|\lesssim_i r^{-1-i}$,
the product rule and \eqref{EQ44} yield
\begin{align}
\begin{split}
 |D^j\Re F|
 &\le |D^jF|
 \lesssim_A |c_ae^{z\beta+h}|
 \sum_{i=0}^j r^{-1-i}|\beta|^{j-i}\\
 &=|F|\sum_{i=0}^j r^{-i}|\beta|^{j-i}
 \lesssim_A |F|\,|\beta|^j,
 \qquad 0\le j\le3.
 \label{EQ50}
 \end{split}
\end{align}
For the Laplacian, \eqref{EQ45} gives
$\Delta\Re F=4\Re\big(\partial(Fz\dbar\beta)\big)$. The factor
$z\dbar\beta$ carries a derivative of the cutoff. In fact, \eqref{EQ49}
implies
$$
 |D^i(z\dbar\beta)|\lesssim r|D^{i+1}\beta|+i|D^i\beta|
 \lesssim_A\frac{k}{\log k}\,r^{-i},\qquad 0\le i\le3.
$$
Applying \eqref{EQ50} to the other factor, we obtain for $0\le j\le2$
\begin{align}
 \begin{split}
 |D^j\Delta\Re F|
 &\lesssim\sum_{i=0}^{j+1}|D^{j+1-i}F|\,|D^i(z\dbar\beta)|\\
 &\lesssim_A\frac{k|F|}{\log k}\sum_{i=0}^{j+1}|\beta|^{j+1-i}r^{-i}
 \lesssim_A\frac{k|F|}{\log k}|\beta|^{j+1}.
 \end{split}\label{EQ51}
\end{align}
For the denominator in \eqref{EQ48}, direct differentiation gives
$$
 \nabla(|\nabla\Re F|^{-2})
 =-\frac{2(D^2\Re F)\nabla\Re F}{|\nabla\Re F|^4}.
$$
Differentiating once more, we obtain
\begin{align*}
 \partial_{ij}(|\nabla\Re F|^{-2})
 &=-2\frac{\sum_{\ell=1}^2
 \big(\partial_{\ell i}\Re F\,\partial_{\ell j}\Re F
 +\partial_\ell\Re F\,\partial_{\ell ij}\Re F\big)}
 {|\nabla\Re F|^4}\\
 &\quad+8\frac{
 \big(\sum_{\ell=1}^2\partial_\ell\Re F\,\partial_{\ell i}\Re F\big)
 \big(\sum_{\ell=1}^2\partial_\ell\Re F\,\partial_{\ell j}\Re F\big)}
 {|\nabla\Re F|^6}.
\end{align*}
By the Cauchy-Schwarz inequality,
\begin{align*}
 |D^2(|\nabla\Re F|^{-2})|
 &\lesssim
 \frac{|D^2\Re F|^2+|\nabla\Re F||D^3\Re F|}
 {|\nabla\Re F|^4}
 +\frac{|\nabla\Re F|^2|D^2\Re F|^2}
 {|\nabla\Re F|^6}\\
 &\lesssim
 \frac{|D^3\Re F|}{|\nabla\Re F|^3}
 +\frac{|D^2\Re F|^2}{|\nabla\Re F|^4}.
\end{align*}
Substitution of \eqref{EQ46} and \eqref{EQ50} yields
$|D^i(|\nabla\Re F|^{-2})|\lesssim_A |F|^{-2}|\beta|^{i-2}$ for
$0\le i\le2$. Multiplying by the derivatives in \eqref{EQ51} gives,
for $0\le j\le2$,
\begin{align}
 \left|D^j\left(\frac{\Delta\Re F}{|\nabla\Re F|^2}\right)\right|
 \lesssim_A\sum_{i=0}^j
 \frac{k|F||\beta|^{i+1}}{\log k}\frac{|\beta|^{j-i-2}}{|F|^2}
 \lesssim_A\frac{k|\beta|^{j-1}}{|F|\log k}.\label{EQ52}
\end{align}
Another application of the product rule, this time to \eqref{EQ47},
gives by \eqref{EQ50} and \eqref{EQ52}
$$
 |D^j g|\lesssim_A\sum_{i=0}^j |F|\,|\beta|^i
 \frac{k|\beta|^{j-i-1}}{|F|\log k}
 \lesssim_A\frac{k}{\log k}|\beta|^{j-1},\qquad j=1,2.
$$
The three terms in \eqref{EQ48} therefore satisfy
\begin{align}
 |V|\lesssim_A \frac{k|\beta|}{\log k}
 +\left(\frac{k}{\log k}\right)^2
 +|F|\,|\beta|\frac{k}{|F|\log k}
 \lesssim_A\frac{k^2}{\log k},\label{EQ53}
\end{align}
where we used $|\beta|\le Ak$ in the last step.

\smallskip\noindent\textbf{(iv) Outer annulus.}
If $b\ne b_*$, the cutoff in \eqref{EQ43} equals one on $d/2<r\le d$.
There $F=R$ is holomorphic, so $\Delta\Re F=0$ and \eqref{EQ47} and
\eqref{EQ48} give $g=V=0$. Thus in both cases $(u,V)=(\Re R,0)$ on
this collar. Extending by the same pair across $\partial B_d$ gives
$u\in C^\infty(\overline B_d)$ and $\supp V\subset\overline B_{d/2}\Subset B_d$.
The estimate~\eqref{EQ42} gives both bounds in the core. Outside the core, $V=0$ if $b=b_*$, while for $b\ne b_*$ the only additional contribution comes from the correction annulus and satisfies~\eqref{EQ53}.
\end{proof}

\section{The boundary and normalization}\label{EQ54}
We start with the proof of the first main theorem.

\begin{proof}[Proof of Theorem~\ref{EQ04}]
By \eqref{EQ29}, Lemma~\ref{EQ31} applies with a
single absolute constant $A$ at every pole. Applying the lemma
in the disjoint disks $B_d(a)$ and setting
$(u,V)=(\Re R,0)$ elsewhere, we obtain a smooth pair
on $\R^2$ satisfying $\Delta u=Vu$, with
$$
 \supp V\subset\{1/4<|z|<1/2\}
 \andand
 c\frac{k^2}{\log k}\le\|V\|_\infty
 \le C\frac{k^2}{\log k}.
$$
Using the cutoff $\eta$ from the proof of Lemma~\ref{EQ31}, set $\zeta(z)=1-\eta(30|z|-24)$. Then $\zeta\in C_c^\infty(\R^2)$, $0\le\zeta\le1$, and
$$
 \zeta=1\quad\text{for }|z|\le4/5
 \andand
 \zeta=0\quad\text{for }|z|\ge5/6.
$$
Set
$\widetilde u=(1-\zeta)+\zeta u$. Since $u=\Re R$
for $|z|\ge1/2$, \eqref{EQ22} gives
$$
 \widetilde u\ge1-|R-1|\ge3/4
    \comma 3/4\le|z|\le7/8
    .
$$
Define
\begin{align}
 V_k(z)=
 \begin{cases}
 V(z),& |z|<3/4,\\[2pt]
 \dfrac{(\Delta\zeta)(\Re R-1)
       +2\nabla\zeta\cdot\nabla\Re R}
       {(1-\zeta)+\zeta\Re R},
       &3/4\le|z|<7/8,\\[7pt]
 0,& |z|\ge7/8.
 \end{cases}\label{EQ55}
\end{align}
Since $\Re R$ is harmonic on the transition annulus, the product rule gives $\Delta\widetilde u=V_k\widetilde u$. The numerator in \eqref{EQ55} vanishes on collars of $|z|=3/4$ and $|z|=7/8$, where the adjacent definitions also give zero. Hence, $V_k\in C_c^\infty(\D)$, with
\begin{align}
 \supp V_k\subset
 \{1/4<|z|<1/2\}\cup\{3/4<|z|<7/8\},
 \label{EQ56}
\end{align}
and $|V_k|\lesssim k(2/3)^N$ on the transition annulus by \eqref{EQ22}, \eqref{EQ23}, and $|\nabla\Re R|=|R'|$. Since $V_k=V$ on the pole disks, $\|V_k\|_\infty\sim k^2/\log k$.

For $|z|\ge1/2$, the product \eqref{EQ18} gives
$$
 |u(z)|\le|R(z)|
 \le\prod_{j=1}^N
 \frac{1}{1-(2\rho_j)^{n_j}}
 =R(1/2).
$$
Let $M_k=\sup_{B_{1/2}}|u|$. By continuity, $M_k\ge u(1/2)=R(1/2)>1$. The preceding bound on $|u|$ outside $B_{1/2}$ then gives $\|u\|_\infty=M_k$. Since $\widetilde u=u$ on $B_{1/2}$ and $0\le\zeta\le1$,
$$
 M_k=\sup_{B_{1/2}}|\widetilde u|\le\|\widetilde u\|_\infty\le\max\{1,\|u\|_\infty\}=M_k.
$$
Set $u_k=M_k^{-1}\widetilde u$.
Then
$$
 \Delta u_k=V_ku_k
 \andand
 \sup_{B_{1/2}}|u_k|=\|u_k\|_\infty=1,
$$
while $u_k=M_k^{-1}>0$ for $|z|\ge7/8$.
Finally, $u_k=M_k^{-1}\Re R$ on~$B_{1/4}$.
By \eqref{EQ21},
$$
 \sup_{B_r}|u_k|
 \le M_k^{-1}(4r)^k\le(4r)^k,
 \qquad 0<r\le1/4.
$$
Furthermore, \eqref{EQ20} gives
$$
 u_k(z)=M_k^{-1} c_k\Re(z^k)
 +O_k(|z|^{k+m}),
 \qquad z\to0.
$$
Along the positive real axis, $u_k(t)/t^k\to M_k^{-1}c_k\ne0$ as $t\downarrow0$. In particular, $\sup_{B_r}|u_k|\ge|u_k(r/2)|\gtrsim_k r^k$ for all sufficiently small $r>0$. Together with the upper bound, this proves $\ord_0u_k=k$.
\end{proof}

\section{Transplantation and sharpness}\label{EQ57}
Since $u_k$ is constant near $\partial D$, we can extend it from an isothermal chart to all of~$X$.
\begin{proof}[Proof of Theorem~\ref{EQ09}]
By compactness, choose finitely many isothermal charts
$\Psi_j\colon B_3\to X$ such that $X=\bigcup_j\Psi_j(B_1)$ and
$$
 \Psi_j^*g=e^{2\phi_j}|\text{d}z|^2
 \andand K^{-1}\le e^{\phi_j}\le K
 \quad\text{on }\overline B_2
$$
for some $K\ge1$. Given $p\in X$, choose $j$ and $y\in B_1$
such that $p=\Psi_j(y)$, and set
$\Psi(z)=\Psi_j(y+z)$ and $\phi(z)=\phi_j(y+z)$.
On $\Psi(\D)$ define
\begin{align}
 U_k\circ\Psi=u_k
 \andand Q_k\circ\Psi=-e^{-2\phi}V_k.
 \label{EQ58}
\end{align}
For $|z|\ge7/8$, we have $u_k=M_k^{-1}$ and, by \eqref{EQ56}, $V_k=0$. Thus the pair in \eqref{EQ58} extends smoothly by $(M_k^{-1},0)$ outside~$\Psi(\D)$. In particular,
$\supp Q_k=\Psi(\supp V_k)\Subset\Psi(\D)$.
The conformal expression of the Laplacian gives
$$
 (\Delta_g U_k+Q_kU_k)\circ\Psi
 =e^{-2\phi}(\Delta u_k-V_ku_k)=0.
$$
Thus $\Delta_g U_k+Q_kU_k=0$ on $X$, while
\eqref{EQ06} yields
$$
 \|Q_k\|_\infty
 \le K^2\|V_k\|_\infty
 \le C_g\frac{k^2}{\log k}
 \andand \|U_k\|_\infty=1.
$$
The metric bounds $K^{-1}|\text{d}z|\le \text{d}s_g\le K|\text{d}z|$
show that every path from $p$ leaving $\Psi(\D)$
has length at least $K^{-1}$, while radial paths
have length at most~$K|z|$. Hence,
$$
 \Psi(B_{r/K})\subset B_r(p)\subset\Psi(B_{Kr}),
 \qquad 0<r<K^{-1},
$$
and consequently
$$
 \sup_{B_{r/K}}|u_k|
 \le\sup_{B_r(p)}|U_k|
 \le\sup_{B_{Kr}}|u_k|.
$$
Together with \eqref{EQ05}, this gives
$\ord_p U_k=k$.
To prove \eqref{EQ11}, set $r_g=(4K)^{-2}$. For $0<r\le r_g$,
we have $Kr\le1/4$ and $4Kr\le\sqrt r$.
Applying \eqref{EQ07} at radius $Kr$ gives
\begin{equation*}
 \sup_{B_r(p)}|U_k|
 \le\sup_{B_{Kr}}|u_k|
 \le(4Kr)^k
 \le r^{k/2}.\qedhere
\end{equation*}
\end{proof}

\begin{proof}[Proof of Corollary~\ref{EQ12}]
For $M>2$ and $0<r<1/2$, the estimate of \cite[Corollary~2.5]{LMNN} gives $\sup_{B_r(p)}|U|\ge r^{C_g\sqrt{M\log M}}\|U\|_\infty$. Hence, $\ord_p U\le C_g\sqrt{M\log M}$.
For the lower bound, choose $\varepsilon>0$ so that $4C_g\varepsilon^2\le1$, with $C_g$ as in \eqref{EQ10}, and take $k=\lfloor\varepsilon\sqrt{M\log M}\rfloor$. For all sufficiently large $M$, we have $k\ge k_g$ and $k\ge\tfrac12\varepsilon\sqrt{M\log M}\ge M^{1/4}$. Thus $\log k\ge\tfrac14\log M$, and Theorem~\ref{EQ09} gives
$$
 \|Q_k\|_\infty\le C_g\frac{k^2}{\log k}\le 4C_g\varepsilon^2 M\le M.
$$
Since $\ord_p U_k=k$, this proves the lower bound in \eqref{EQ13} with $c_g=\varepsilon/2$.
\end{proof}

\section{Decay at infinity}\label{EQ59}

The decay in \eqref{EQ16} comes from the factors whose poles lie on the inner circles.

\begin{proof}[Proof of Theorem~\ref{EQ15}]
We repeat the construction of Theorem~\ref{EQ04} on infinitely many circles. The $j$th circle has radius $r_j=j^2/\sqrt{\log j}$ and carries $j$ equally spaced poles, and each pole is regularized by Lemma~\ref{EQ31} with $k=j^2$ in the rescaled variable $w=z/r_j$. With this choice of $r_j$, the potential stays bounded after rescaling, while the solution has size $e^{-cj^2}$ near the $j$th circle, where $j^2$ is comparable to $r_j\sqrt{\log r_j}$. The proof has five steps.
\begin{enumerate}
\item[(i)] Construct a convergent product $P$ with $j$ equally spaced poles on the $j$th circle.
\item[(ii)] Separate the pole disks and verify the bounds on $b,h$ in Lemma~\ref{EQ31} with $k=j^2$.
\item[(iii)] Patch the regularized solutions to $\Re P$ and bound the potential.
\item[(iv)] Extract $e^{-cj^2}$ from the inner factors and estimate the solution outside and inside the pole disks.
\item[(v)] Translate a maximum point of $|u|$ to the origin and normalize the solution and potential.
\end{enumerate}

\smallskip\noindent\textbf{(i) The infinite product.}
For an integer $J\ge3$, to be fixed at the end of (ii), set
\begin{align}
 r_j=\frac{j^2}{\sqrt{\log j}}\quad(j\ge2)
 \andand P(z)=\prod_{i=J}^\infty\frac{1}{1-(z/r_i)^i}.\label{EQ60}
\end{align}
The logarithmic derivative of $\tau^2/\sqrt{\log\tau}$ is $(2-1/2\log\tau)/\tau$, which lies between $1/\tau$ and $2/\tau$ for $\tau\ge2$. Integrating from $i$ to $j$, we get
\begin{align}
 \log\frac ji\le\log\frac{r_j}{r_i}
 =\int_i^j\left(2-\frac1{2\log\tau}\right)\frac{\text{d}\tau}{\tau}
 \le2\log\frac ji\comma 2\le i<j.\label{EQ61}
\end{align}
In particular, the radii are strictly increasing and tend to infinity.

To prove that the product converges, let $R>0$ and choose $\tilde J=\tilde J(R)\ge J$ so that $r_i\ge2R$ for $i\ge\tilde J$. Since $|\log(1-t)|\le|t|/(1-|t|)$ for $|t|<1$, we have
\begin{align}
 \sum_{i=\tilde J}^\infty\sup_{|z|\le R}\big|\log(1-(z/r_i)^i)\big|
 \le\sum_{i=\tilde J}^\infty\frac{(R/r_i)^i}{1-(R/r_i)^i}
 \le2\sum_{i=\tilde J}^\infty 2^{-i}<\infty.\label{EQ62}
\end{align}
Hence, $\prod_{i\ge\tilde J}(1-(z/r_i)^i)$ converges uniformly on $B_R$ to a holomorphic function without zeros. Since $R$ is arbitrary, the denominator $\prod_{i\ge J}(1-(z/r_i)^i)$ of $P$ is an entire function, and its zeros are the points
$$
 r_j\alpha\comma\text{where}\quad j\ge J,\quad \alpha^j=1.
$$
Each of them is a simple zero of the $j$th factor and of no other factor, because the radii are distinct. Therefore, $P$ is meromorphic on $\C$ with simple poles exactly at these points. It has no zeros, and $P(0)=1$.

\smallskip\noindent\textbf{(ii) The pole disks.}
All constants in this step are absolute. In particular, they do not depend on $J$.

\smallskip\noindent\textit{Separation of the circles.}
We first compare the radii. By \eqref{EQ61} and the mean value theorem for the logarithm,
\begin{align}
 \frac i2\left|\log\frac{r_j}{r_i}\right|
 \ge\frac i2\left|\log\frac ji\right|
 \ge\frac{i|j-i|}{2\max\{i,j\}}\ge\frac{i|j-i|}{2(i+|j-i|)}
 \ge\frac14\min\{i,|i-j|\}\comma i\ne j.\label{EQ63}
\end{align}
For $j\ge J$, consider the annulus
$$
 A_j=\big\{z:\sqrt{r_{j-1}r_j}\le|z|\le\sqrt{r_jr_{j+1}}\big\}.
$$
The annuli $A_j$ have pairwise disjoint interiors. Fix $j\ge J$ and $z\in A_j$, and for $i\ge J$ with $i\ne j$ write
$$
 t_i=(z/r_i)^i \andand m_i=\min\{i,|i-j|\}.
$$
Since the radii increase, $|z|^2\ge r_ir_j$ for $i<j$ and $|z|^2\le r_ir_j$ for $i>j$. Together with \eqref{EQ63}, this gives
\begin{align}
 \begin{split}
 |t_i|&\ge\left(\frac{r_j}{r_i}\right)^{i/2}\ge e^{m_i/4}\comma J\le i<j,\\
 |t_i|&\le\left(\frac{r_j}{r_i}\right)^{i/2}\le e^{-m_i/4}\comma i>j.
 \end{split}\label{EQ64}
\end{align}
These are the analogs of the bounds on $t=(z/\rho_i)^{n_i}$ in Section~\ref{EQ17}. To sum them over $i$, note that for each integer $m\ge1$ there are at most three indices $i$ with $m_i=m$, namely $m$, $j-m$, and $j+m$, and each of them is at most $j+m$.

\smallskip\noindent\textit{Separation of the disks.}
Fix $j\ge J$ and $\alpha$ with $\alpha^j=1$, and write $w=z/r_j$. For $k=j^2$, the radius in Lemma~\ref{EQ31} is $d_j=1/[64(j+1)]$. If $|w-\alpha|\le8d_j$, then $\big||w|-1\big|\le8d_j<1/8j$. On the other hand, \eqref{EQ61} gives
$$
 \sqrt{\frac{r_{j+1}}{r_j}}\ge\sqrt{1+\frac1j}\ge1+\frac1{3j}
 \andand
 \sqrt{\frac{r_{j-1}}{r_j}}\le\sqrt{1-\frac1j}\le1-\frac1{2j}.
$$
Hence, the closed disk $\overline{B_{8r_jd_j}(r_j\alpha)}$ lies in the interior of~$A_j$, and disks centered on different circles are disjoint. On the same circle, two distinct centers satisfy $|\alpha-\alpha'|\ge2\sin(\pi/j)\ge4/j>16d_j$ in the $w$ variable. Thus the disks $B_{8r_jd_j}(r_j\alpha)$ are pairwise disjoint.

\smallskip\noindent\textit{The local form of $P$.}
Omitting the $j$th factor, write
\begin{align}
 H_j(w)=\prod_{\substack{i\ge J\\i\ne j}}\frac1{1-(r_jw/r_i)^i}
 \andand P(r_jw)=\frac{H_j(w)}{1-w^j}.\label{EQ65}
\end{align}
As in (i), $H_j$ is meromorphic without zeros, and its poles lie on the circles $|w|=r_i/r_j$ with $i\ne j$. By \eqref{EQ64}, none of these circles meets $r_j^{-1}A_j$. Thus $H_j$ is holomorphic and nonvanishing on a neighborhood of $\overline{B_{8d_j}(\alpha)}$. Since $1-w^j$ has a simple zero at $\alpha$ with derivative $-j\alpha^{j-1}$, the residue of $P(r_jw)$ at $w=\alpha$ is
\begin{align}
 c_\alpha=-\frac{H_j(\alpha)}{j\alpha^{j-1}}=-\frac{\alpha}{j}H_j(\alpha)\ne0.\label{EQ66}
\end{align}
As in Section~\ref{EQ17}, the function $(w-\alpha)P(r_jw)/c_\alpha$ is holomorphic and nonvanishing on $B_{8d_j}(\alpha)$ and equals one at~$\alpha$. Let $p$ be its holomorphic logarithm with $p(\alpha)=0$. By \eqref{EQ65},
\begin{align}
 p'(w)=\frac{H_j'(w)}{H_j(w)}
 +\left(\frac1{w-\alpha}+\frac{jw^{j-1}}{1-w^j}\right)
 \comma w\ne\alpha.\label{EQ67}
\end{align}
With $\sigma_j=\sum_{i=J}^{j-1}i$, we claim, in analogy with \eqref{EQ24}, that
\begin{align}
 p'(w)=-\frac{\sigma_j}{w}+O(j)\comma |w-\alpha|\le4d_j.\label{EQ68}
\end{align}

\smallskip\noindent\textit{The factors with $i\ne j$.}
Let $|w-\alpha|\le8d_j$, so that $z=r_jw\in A_j$, and let $t_i=(r_jw/r_i)^i$ be as in \eqref{EQ64}. By \eqref{EQ62}, the product for $1/H_j$ converges locally uniformly, so $H_j$ may be differentiated logarithmically term by term. The $i$th factor contributes
$$
 w\frac{\text{d}}{\text{d}w}\log\frac1{1-t_i}
 =\frac{it_i}{1-t_i}
 =-i+\frac{i}{1-t_i}.
$$
Using the first expression for $i>j$ and the second for $i<j$, we obtain
\begin{align*}
 w\frac{H_j'(w)}{H_j(w)}+\sigma_j
 =\sum_{i=J}^{j-1}\frac{i}{1-t_i}+\sum_{i=j+1}^\infty\frac{it_i}{1-t_i}.
\end{align*}
By \eqref{EQ64},
\begin{align*}
 \left|\frac1{1-t_i}\right|&\le\frac1{|t_i|-1}\le\frac1{e^{m_i/4}-1}
 \comma i<j,\\
 \left|\frac{t_i}{1-t_i}\right|&\le\frac{|t_i|}{1-|t_i|}\le\frac1{e^{m_i/4}-1}
 \comma i>j.
\end{align*}
Summing over $i\ne j$ and using the count following \eqref{EQ64}, we get
\begin{align*}
 \left|w\frac{H_j'(w)}{H_j(w)}+\sigma_j\right|
 \le\sum_{\substack{i\ge J\\i\ne j}}\frac{i}{e^{m_i/4}-1}
 \le3\sum_{m=1}^\infty\frac{j+m}{e^{m/4}-1}\lesssim j.
\end{align*}

\smallskip\noindent\textit{The factor with $i=j$.}
We apply \eqref{EQ28} with $k=j^2$. In Section~\ref{EQ17} we then have $N=j$, $d=d_j$, and $n_1=j$. Thus the first factor of \eqref{EQ18} is $R_1(z)=z^j/(z^j-\rho_1^j)$, and $a=\rho_1\alpha$ is one of its poles. Let $0<|w-\alpha|\le4d_j$ and $z=\rho_1w$. Then $|z-a|=\rho_1|w-\alpha|<4d$ and $R_1(\rho_1w)=w^j/(w^j-1)$, so the chain rule and \eqref{EQ28} give
$$
 \left|\frac1{w-\alpha}+\frac{\text{d}}{\text{d}w}\log\frac{w^j}{w^j-1}\right|
 =\rho_1\left|\frac1{z-a}+\frac{R_1'(z)}{R_1(z)}\right|\lesssim j.
$$
The logarithmic derivative on the left equals $j/w+jw^{j-1}/(1-w^j)$, and $|w|>1/2$. Hence, the term in parentheses in \eqref{EQ67} is~$O(j)$. Together with the bound for the factors with $i\ne j$, this proves \eqref{EQ68} for $w\ne\alpha$, and by continuity also at $w=\alpha$.

\smallskip\noindent\textit{The bounds on $b$ and $h$.}
Set $b=p'(\alpha)$ and $h(w)=p(w)-b(w-\alpha)$. Then $h(\alpha)=h'(\alpha)=0$, and since $\sigma_j\le j^2/2$, \eqref{EQ68} gives $|b|\le\sigma_j+Cj\lesssim j^2$. Applying \eqref{EQ68} at both $w$ and $\alpha$, we obtain
\begin{align*}
 |h'(w)|=|p'(w)-p'(\alpha)|
 \le\sigma_j\left|\frac1w-\frac1\alpha\right|+Cj
 \lesssim j^2d_j+j\lesssim j
 \comma |w-\alpha|\le4d_j.
\end{align*}
Since $h(\alpha)=0$, integration along the segment from $\alpha$ to $w$ yields $\|h\|_{L^\infty(B_{2d_j}(\alpha))}\lesssim jd_j\lesssim1$. Thus, as in \eqref{EQ29},
\begin{align}
 P(r_jw)=\frac{c_\alpha}{w-\alpha}e^{b(w-\alpha)+h(w)}
 \andand \frac{|b|}{j^2}+\|h\|_{L^\infty(B_{2d_j}(\alpha))}\le A.\label{EQ69}
\end{align}
Here, $A\ge1$ is independent of $J,j,\alpha$. Fix $J$ with $J^2\ge k_0(A)$ so that Lemma~\ref{EQ31} applies at every pole.

\smallskip\noindent\textbf{(iii) The global solution.}
By \eqref{EQ69}, Lemma~\ref{EQ31} applies in the variable $w$ with $k=j^2$, $d=d_j$, $a=\alpha$, and $R(w)=P(r_jw)$. Let $(v,W)$ be the pair on $B_{d_j}(\alpha)$ constructed in its proof. Thus $\Delta_wv=Wv$, the pair agrees with $(\Re P(r_jw),0)$ near $\partial B_{d_j}(\alpha)$, and $\|W\|_\infty\le C_Aj^4/\log(j^2)$. We refer to $B_{r_jd_j}(r_j\alpha)$ as the pole disks. As in Theorem~\ref{EQ04}, set
\begin{align}
 (u(z),V(z))=
 \begin{cases}
 (v(w),r_j^{-2}W(w)),&z=r_jw,\ |w-\alpha|<d_j,\\
 (\Re P(z),0),&\text{outside the pole disks}.
 \end{cases}\label{EQ70}
\end{align}
By (ii), the pole disks are pairwise disjoint, and those centered on the $j$th circle lie in~$A_j$. In particular, only finitely many of them meet a given bounded set, and the origin lies outside all of them. Since the two formulas in \eqref{EQ70} agree near the boundary of each pole disk, $u$ and $V$ are smooth and real on $\R^2$, with $u(0)=P(0)=1$. On each pole disk $\Delta u=r_j^{-2}\Delta_wv=Vu$, and $\Re P$ is harmonic elsewhere. Finally, the choice of the radii gives
\begin{align*}
 \|V\|_\infty\le C_A\sup_{j\ge J}\frac{j^4}{r_j^2\log(j^2)}=\frac{C_A}{2}.
\end{align*}

\smallskip\noindent\textbf{(iv) Decay on the annuli.}
Throughout this step $j\ge4J$, and the constants do not depend on~$j$. By (ii), the only pole disks that meet $A_j$ are those centered on the $j$th circle. We show that, with $c=(\log2)/64$,
\begin{align}
 \sup_{A_j}|u|\lesssim e^{-cj^2/2}.\label{EQ71}
\end{align}

\smallskip\noindent\textit{The inner factors.}
Let $z=r_jw\in A_j$, and let $t_i$ and $m_i$ be as in \eqref{EQ64}. Since $|1-t_i|\ge|t_i|(1-|t_i|^{-1})$ for $i<j$ and $|1-t_i|\ge1-|t_i|$ for $i>j$, \eqref{EQ65}, \eqref{EQ64}, and the count following \eqref{EQ64} give
\begin{align*}
 |H_j(w)|
 \le\prod_{i=J}^{j-1}\frac1{|t_i|}
 \prod_{\substack{i\ge J\\i\ne j}}\frac1{1-e^{-m_i/4}}
 \le\prod_{m=1}^\infty\frac1{(1-e^{-m/4})^3}\,
 \prod_{i=J}^{j-1}\frac1{|t_i|}.
\end{align*}
The product over $m$ is finite since $\sum_{m\ge1}e^{-m/4}<\infty$. In the last product, every factor is at most one. Moreover, for $j/4\le i\le j/2$ we have $i\ge J$, and \eqref{EQ64} and \eqref{EQ61} give $|t_i|\ge(r_j/r_i)^{i/2}\ge(j/i)^{i/2}\ge2^{i/2}\ge2^{j/8}$. Since there are at least $j/4-1\ge j/8$ such indices,
\begin{align}
 |H_j(w)|\lesssim 2^{-j^2/64}=e^{-cj^2}
 \comma z=r_jw\in A_j.\label{EQ72}
\end{align}

\smallskip\noindent\textit{Outside the pole disks.}
Let $z=r_jw\in A_j$ lie outside the pole disks. Then $u(z)=\Re P(z)$, and $|w-\alpha|\ge d_j$ whenever $\alpha^j=1$, so that $|1-w^j|=\prod_{\alpha^j=1}|w-\alpha|\ge d_j^j$. Since $j\log(1/d_j)=O(j\log j)$, \eqref{EQ65} and \eqref{EQ72} give
\begin{align*}
 |u(z)|\le|P(z)|=\frac{|H_j(w)|}{|1-w^j|}
 \lesssim d_j^{-j}e^{-cj^2}\lesssim e^{-cj^2/2}.
\end{align*}

\smallskip\noindent\textit{Inside the pole disks.}
Consider the pole disk centered at $r_j\alpha$, and let $k=j^2$. Since $r_j\alpha\in A_j$, \eqref{EQ66} and \eqref{EQ72} give $|c_\alpha|\lesssim j^{-1}e^{-cj^2}$. We claim that
\begin{align}
 \sup_{B_{d_j}(\alpha)}|v|\lesssim\frac{k|c_\alpha|}{\sqrt{\log k}}\,e^{|b|d_j}.\label{EQ73}
\end{align}
To prove \eqref{EQ73}, we return to the proof of Lemma~\ref{EQ31} and write $r=|w-\alpha|$. In the core $r<\sqrt{\log k}/k$, \eqref{EQ36} gives $|v|\le k|c_\alpha|(\log k)^{-1/2}\sup_{t\ge0}tq(t)$, and the supremum is finite since, by \eqref{EQ35}, $q$ is bounded and $q(t)=t^{-2}$ for $t\ge1/4$. For $\sqrt{\log k}/k\le r\le d_j$, we have $v=e^g\Re F$, where $F$ is given by \eqref{EQ32} or \eqref{EQ43}, and $g=0$ outside the correction annulus. The linear coefficient in the exponent of $F$ is $b_*$, $\beta(r)$, or $b$, and thus has modulus at most $|b|$. Since $|h|\le A$ by \eqref{EQ69}, we get
\begin{align*}
 |F|\le\frac{|c_\alpha|}{r}e^{|b|r+A}
 \le\frac{k|c_\alpha|}{\sqrt{\log k}}e^{|b|d_j+A}.
\end{align*}
On the correction annulus, \eqref{EQ47}, \eqref{EQ52} at order zero, and \eqref{EQ44} give
$$
 |g|\le\frac{|F|}{2}\left|\frac{\Delta\Re F}{|\nabla\Re F|^2}\right|
 \lesssim\frac{k}{|\beta|\log k}\lesssim\frac1{\sqrt{\log k}},
$$
so $e^g\lesssim1$, and \eqref{EQ73} follows. Since $k=j^2$ and $|b|d_j\le Aj^2d_j\le Aj$, \eqref{EQ73} gives
\begin{align*}
 \sup_{B_{d_j}(\alpha)}|v|
 \lesssim j e^{-cj^2+Aj}\lesssim e^{-cj^2/2}.
\end{align*}
By \eqref{EQ70}, the bounds outside and inside the pole disks prove \eqref{EQ71}.

\smallskip\noindent\textit{Decay in $|z|$.}
The annuli $A_j$ with $j\ge4J$ cover all $z$ with $|z|\ge\sqrt{r_{4J-1}r_{4J}}$. If $z\in A_j$, then $|z|\le r_{j+1}\le(j+1)^2$, and since $t\sqrt{\log t}$ increases for $t\ge1$,
$$
 |z|\sqrt{\log|z|}\le r_{j+1}\sqrt{\log r_{j+1}}
 \le\frac{(j+1)^2}{\sqrt{\log(j+1)}}\sqrt{2\log(j+1)}\le3j^2.
$$
Thus \eqref{EQ71} gives $|u(z)|\lesssim e^{-(c/6)|z|\sqrt{\log|z|}}$ on $A_j$. After decreasing $c$, we obtain
\begin{align}
 |u(z)|\le e^{-c|z|\sqrt{\log|z|}}\qquad\text{for all sufficiently large }|z|.\label{EQ74}
\end{align}

\smallskip\noindent\textbf{(v) Normalization.}
Since $u$ is continuous, $u(0)=1$, and $u$ tends to zero at infinity by \eqref{EQ74}, we have $|u(p)|=\|u\|_\infty\ge1$ for some $p\in\R^2$. With $\lambda=\max\{1,\sqrt{\|V\|_\infty}\}$ set
\begin{align*}
 \widetilde u(z)=\frac{u(p+z/\lambda)}{u(p)}\andand
 \widetilde V(z)=\lambda^{-2}V(p+z/\lambda).
\end{align*}
Then $\Delta\widetilde u=\widetilde V\widetilde u$, $\widetilde u(0)=\|\widetilde u\|_\infty=1$, and $\|\widetilde V\|_\infty\le1$. Moreover, $|\widetilde u(z)|\le|u(p+z/\lambda)|$, while $|p+z/\lambda|=|z|/\lambda+O(1)$ and $\log|p+z/\lambda|=\log|z|+O(1)$ as $|z|\to\infty$. Hence, \eqref{EQ74} gives \eqref{EQ16} for $\widetilde u$ with $c/2\lambda$ in place of~$c$.
\end{proof}

\section*{Acknowledgments}
The authors have used ChatGPT and Claude to assist with the proofs and writing. The authors take full responsibility for the content.
The authors were supported in part 
by the Simons grant SFI-MPS-TSM-00025796.


\begin{thebibliography}{99}

\bibitem{Aronszajn}
N.~Aronszajn, \emph{A unique continuation theorem for solutions of elliptic partial differential equations or inequalities of second order}, J. Math. Pures Appl. (9) \textbf{36} (1957), 235--249.

\bibitem{BakriC1}
L.~Bakri, \href{https://doi.org/10.1512/iumj.2012.61.4713}\emph{Quantitative uniqueness for Schr\"odinger operator}, Indiana Univ. Math. J. \textbf{61} (2012), 1565--1580.

\bibitem{BakriComplex}
L.~Bakri, \href{https://doi.org/10.1080/03605302.2012.736912}\emph{Carleman estimates for the Schr\"odinger operator. Applications to quantitative uniqueness}, Comm. Partial Differential Equations \textbf{38} (2013), 69--91.

\bibitem{BK}
J.~Bourgain and C.~E.~Kenig, \href{https://doi.org/10.1007/s00222-004-0435-7}\emph{On localization in the continuous Anderson-Bernoulli model in higher dimension}, Invent. Math. \textbf{161} (2005), 389--426.

\bibitem{CKW1}
G.~Camliyurt, I.~Kukavica, and F.~Wang, \href{https://doi.org/10.1007/s00209-018-2081-6}\emph{On quantitative uniqueness for elliptic equations}, Math. Z. \textbf{291} (2019), 227--244.

\bibitem{CKW2}
G. Camliyurt, I. Kukavica, and F. Wang,
On localization and quantitative uniqueness for elliptic
partial differential equations,
in \emph{Partial Differential Equations in Fluid Mechanics},
C. L. Fefferman, J. C. Robinson, and J. L. Rodrigo (eds.),
London Math. Soc. Lecture Note Ser., vol.~452,
Cambridge Univ. Press, Cambridge, 2018, 68--96.

\bibitem{Carleman}
T.~Carleman, \emph{Sur un probl\`eme d'unicit\'e pour les syst\`emes d'\'equations aux d\'eriv\'ees partielles \`a deux variables ind\'ependantes}, Ark. Mat. Astr. Fys. \textbf{26B} (1939), no.~17, 1--9.

\bibitem{Davey}
B.~Davey, \href{https://doi.org/10.1016/j.jfa.2020.108566}\emph{Quantitative unique continuation for Schr\"odinger operators}, J. Funct. Anal. \textbf{279} (2020), no.~4, 108566.

\bibitem{DKW}
B.~Davey, C.~E.~Kenig, and J.-N.~Wang, \emph{On Landis' conjecture in the plane when the potential has an exponentially decaying negative part}, Algebra i Analiz \textbf{31} (2019), no.~2, 204--226.

\bibitem{DZ}
B.~Davey and J.~Zhu, \href{https://doi.org/10.1080/03605302.2019.1629957}\emph{Quantitative uniqueness of solutions to second-order elliptic equations with singular lower order terms}, Comm. Partial Differential Equations \textbf{44} (2019), 1217--1251.

\bibitem{DF}
H.~Donnelly and C.~Fefferman, \href{https://doi.org/10.1007/BF01393691}\emph{Nodal sets of eigenfunctions on Riemannian manifolds}, Invent. Math. \textbf{93} (1988), 161--183.

\bibitem{FGJL}
S.~Feng, Z.~Gan, R.~Jiang, and F.~Lin, \emph{Vanishing orders of Dirichlet solutions to the Schr\"odinger equation in dimensions three and higher}, \href{https://arxiv.org/abs/2610.01615}{arXiv:2610.01615} (2026).

\bibitem{FSR}
A.~Fern\'andez-Bertolin, L.~Roncal, and D.~Stan,
\href{https://doi.org/10.17983/301255}\emph{Landis' conjecture: a survey},
in \emph{Harmonic Analysis and Nonlinear Partial Differential Equations},
RIMS K\^oky\^uroku Bessatsu \textbf{B99} (2026), 35--78.



\bibitem{FK}
N.~D.~Filonov and S.~T.~Krymskii, \href{https://doi.org/10.1134/S1061920824040058}\emph{On the Landis conjecture in a cylinder}, Russ. J. Math. Phys. \textbf{31} (2024), 645--665.

\bibitem{FI}
R.~L.~Frank and P.~Ivanisvili, \emph{Counterexamples to the Landis conjecture in dimensions three and higher}, \href{https://arxiv.org/abs/2608.00802}{arXiv:2608.00802} (2026).

\bibitem{GL}
N.~Garofalo and F.-H.~Lin, \href{https://doi.org/10.1512/iumj.1986.35.35015}\emph{Monotonicity properties of variational integrals, $A_p$ weights and unique continuation}, Indiana Univ. Math. J. \textbf{35} (1986), 245--268.

\bibitem{JK}
D.~Jerison and C.~E.~Kenig, \href{https://doi.org/10.2307/1971205}\emph{Unique continuation and absence of positive eigenvalues for Schr\"odinger operators}, Ann. of Math. (2) \textbf{121} (1985), 463--488.

\bibitem{Kenig}
C.~E.~Kenig, \href{https://media.accademiaxl.it/memorie/S5-VXXIX-P1-2-2005/Kenig231-242.pdf}\emph{Some recent quantitative unique continuation theorems}, Rend. Accad. Naz. Sci. XL Mem. Mat. Appl. \textbf{29} (2005), no.~1, 231--242. Also in \href{https://www.numdam.org/item/SEDP_2005-2006____A20_0/}\emph{S\'eminaire \'Equations aux d\'eriv\'ees partielles}, \'Ecole Polytechnique (2005--2006), Exp.~no.~XX, 1--10.

\bibitem{KSW}
C.~E.~Kenig, L.~Silvestre, and J.-N.~Wang, \href{https://doi.org/10.1080/03605302.2014.978015}\emph{On Landis' conjecture in the plane}, Comm. Partial Differential Equations \textbf{40} (2015), 766--789.

\bibitem{KL}
V.~A.~Kondrat'ev and E.~M.~Landis, \href{https://doi.org/10.1007/978-3-642-58173-1_2}\emph{Qualitative theory of second order linear partial differential equations}, in \emph{Partial Differential Equations III}, Encyclopaedia Math. Sci., vol.~32, Springer, Berlin, 1991, 87--192.

\bibitem{Kukavica}
I.~Kukavica, \href{https://doi.org/10.1215/S0012-7094-98-09111-6}\emph{Quantitative uniqueness for second-order elliptic operators}, Duke Math. J. \textbf{91} (1998), no.~2, 225--240.

\bibitem{LBS}
K.~Le Balc'h and D.~A.~Souza, \href{https://doi.org/10.2140/apde.2026.19.1061}\emph{Quantitative unique continuation for real-valued solutions to second order elliptic equations in the plane}, Anal. PDE \textbf{19} (2026), 1061--1106.

\bibitem{LMNN}
A.~Logunov, E.~Malinnikova, N.~Nadirashvili, and F.~Nazarov, \href{https://doi.org/10.1007/s00222-025-01340-1}\emph{The Landis conjecture on exponential decay}, Invent. Math. \textbf{241} (2025), 465--508.

\bibitem{Meshkov}
V.~Z.~Meshkov, \href{https://doi.org/10.1070/SM1992v072n02ABEH001414}\emph{On the possible rate of decay at infinity of solutions of second order partial differential equations}, Math. USSR-Sb. \textbf{72} (1992), 343--361.

\bibitem{Zhu}
J.~Zhu, \href{https://doi.org/10.1353/ajm.2016.0027}\emph{Quantitative uniqueness of elliptic equations}, Amer. J. Math. \textbf{138} (2016), 733--762.

\end{thebibliography}
\end{document}